\documentclass[a4paper, 11pt]{amsart}
\usepackage{amssymb,amsmath,txfonts,mathrsfs}
\newtheorem{theorem}{Theorem}[section]

\newtheorem{remark}[theorem]{Remark}

\newtheorem{cor}[theorem]{Corollary}

\numberwithin{equation}{section}
\begin{document}

\title[local pinching estimates in Ricci flow]{Local pinching estimates in $3$-dim Ricci flow}
\author{Bing-Long Chen, Guoyi Xu, Zhuhong Zhang}
\address{Bing-Long Chen\\ Department of Mathematics\\Sun Yat-Sen University\\Guangzhou\\P. R. China}
\email{mcscbl@mail.sysu.edu.cn}
\address{Guoyi Xu\\ Mathematical Sciences Center\\Tsinghua University, Beijing\\100084\\P. R. China}
\email{gyxu@math.tsinghua.edu.cn}
\address{Zhuhong Zhang\\ Department of Mathematics\\South China Normal Univeristy\\Guangzhou\\P. R. China}
\email{juhoncheung@sina.com}

\begin{abstract}
We study curvature pinching estimates of Ricci flow on complete $3$-dimensional manifolds without bounded curvature assumption. We will derive some general curvature conditions which are preserved on any complete solution of $3$-dim Ricci flow, these conditions include nonnegative Ricci curvature and sectional curvature  as  special cases. A local version of Hamilton-Ivey estimates is also obtained.
\\[3mm]
Mathematics Subject Classification: 35K40, 53C44
\end{abstract}

\maketitle

\section{Introduction}

Ricci flow on noncompact complete manifolds has received much attention since W-X Shi's fundamental work in 1990's (see \cite{Shi1}, \cite{Shi2}). Many progress in the field are achieved under bounded curvature assumption, which essentially support the use of maximum principle globally on the whole manifold. Without (local) bounded curvature assumption for the solution of the Ricci flow, very few (localized) estimates have been obtained. Note that Shi's local gradient estimates depends on the local curvature bound of the solution.

 In \cite{Chen}, the first author obtained some local estimates on lower bounds of curvature operators in $3$-dim Ricci flow only assuming the completeness of the solution. The main localization techniques used there is the one  employed by P. Li and S-T Yau in \cite{LY} for the parabolic equations on manifolds. The modification for the Ricci flow situation is to replace the Laplacian comparison theorem by the comparison result of Perelman on the distance function. As a corollary,  it was shown that nonnegative sectional curvature is preserved under $3$ dimensional complete smooth solution to the Ricci flow. For three-manifolds, the nonnegative curvature operator is equivalent to the nonnegative sectional curvature.

In dimension 3, Ricci curvature also determines the whole curvature operator.  One purpose of this note is to generalize the result in \cite{Chen} to the Ricci curvature case.  We will show (see Theorem \ref{thm 1.2})that for any fixed nonnegative constants $a$, $b, c\geq 0$, $a R+b\cdot \min{Rc}+ c\cdot \min{\mathscr M}\geq 0$ is preserved on any $3$-dim complete solution to the Ricci flow without bounded curvature assumption on the solution, where $\mathscr M_{ij}= Rg_{ij}- 2R_{ij}$ is the curvature tensor of $M^3$.  More precisely,
 \begin{theorem}\label{thm 1.1}
{Given a smooth complete solution of Ricci flow $(M^3, g(t))$ on $[0, T]$, let $\lambda\geq \mu\geq \nu$ are eigenvalues of curvature operator matrix $\mathscr M_{ij}$. Suppose for some nonnegative constants $a ,b, c\geq 0,$ we have $[a(\lambda+\mu+\nu)+b(\mu+\nu)+c\nu](g(0))\geq 0$ at time $t=0,$ then $[a(\lambda+\mu+\nu)+b(\mu+\nu)+c\nu](g(t))\geq 0$ for $t>0$.
}
\end{theorem}

From this general curvature pinching result, we find  that $R\geq 0,$ $Rc\geq 0$ and $\mathscr M\geq 0$ are all preserved under any $3$-dim  complete Ricci flow solution. In particular, we have
\begin{theorem}\label{thm 1.2}
{$(M^3, g(t))$ is a smooth complete solution of Ricci flow on $[0, T]$, if $Rc(g(0))\geq 0$, then $Rc(g(t))\geq 0$.
}
\end{theorem}

  Note that  in \cite{Ha82} R. Hamilton proved that $Rc\geq 0$ is preserved for the Ricci flow on closed 3-manifolds.

 In 3-dim Ricci flow, the magic Hamilton-Ivey pinching estimate (or its improvement) states that the least eigenvalues of the curvature operator are not comparable with the largest eigenvalues on the high curvature region. More precisely, if $\nu(x,0)\geq -1$(this is always possible by scaling the initial metric), then for $t>0,$ we have \begin{equation}\label{hip}R\geq (-\nu)(-3+\ln ((1+t)(-v)))\end{equation}  when $\nu<0,$ where $\lambda\geq\mu\geq\nu$ are eigenvalues of the curvature operator, $R=\lambda+\mu+\nu.$ Note that (\ref{hip}) was only proved  (see \cite{Ha99}) previously on compact solutions or complete noncompact solutions with bounded curvature.  The result in \cite{Chen} is just  the linearized and localized version of the above estimate.  The second major purpose of this note is to derive  a genuine local version of the (improved) Hamilton-Ivey pinching estimate (\ref{hip}).
\begin{theorem}\label{thm 1.3}
Let {$(M^3, g(t))$ be a 3-dim complete smooth solution of the Ricci flow on $[0, T]$. For any fixed $r_0,K>0,$ assume $Rc(x, t)\leq 2r_0^{-2}$ for $x\in B_{t}(x_0, r_0)$, $t\in [0, T]$, and $\nu(x, 0)\geq -K$ on $B_{0}(x_0, 2Ar_0)$ at time $t=0$, where $A\geq \frac{C}{3r_0^2}\Big(\frac{1}{K}+ T\Big),$ $C> 0$ is some universal constant.  Then at any point $(x, t)\in B_{t}(x_0, Ar_0)$, $t\in [0, T]$ where $\nu(x, t)< 0$, we have
\begin{equation}\label{1.3.1}
{\frac{R}{-\nu}- \ln (-\nu)+ 3+ \ln\frac{K}{1+ Kt}\geq \min\Big\{-\frac{8640(1+ KT)}{KAe^2r_0^2},\ -\frac{31104(1+ KT)}{Ke^3(Ar_0)^2} \Big\}.
}
\end{equation}
}
\end{theorem}

We remark that in Theorem \ref{thm 1.3}, the completeness of the solution is not necessary, we only need the condition that all balls $B_{t}(x_0,r_0)$ and $B_t(x_0,2Ar_0)$ are compactly contained in the manifold.

\begin{cor}\label{thm 1.4}
Let {$(M^3, g(t))$ be a 3-dim complete smooth solution of the Ricci flow on $[0, T]$. For any fixed $ 0<K\leq \infty$, if $\nu(x, 0)\geq -K$ on $M^3$ at time $t=0.$ Then at any point $(x, t)\in M^3\times [0, T]$ with $\nu(x, t)< 0$, we have
\begin{equation}\label{1.3.2}
{\frac{R}{-\nu}- \ln (-\nu)+ 3+ \ln\frac{K}{1+ Kt}\geq 0.
}
\end{equation}
}
\end{cor}
A special case of  Corollary \ref{thm 1.4}  is

\begin{cor}\label{thm 1.5}
Let $(M^3, g(t))$ be a 3-dim complete smooth solution of the Ricci flow on $[0, T]$.  Then  we have
\begin{equation}\label{1.3.3}
R\geq (-\nu)(-3+\ln (t(-\nu)))
\end{equation}
at any point $(x, t)\in M^3\times [0, T]$ with $\nu(x, t)< 0$.
\end{cor}
In the end of this section, we mention one interesting application of Theorem \ref{thm 1.2}. We may generalize the strong uniqueness theorem of the first author in \cite{Chen} to the case where the initial manifold has only bounded nonnegative Ricci curvature and a uniform injectivity radius lower bound. See Theorem \ref{thm 4.1} for the precise statement.

\section{The preservation of $Rc\geq 0$ and local Hamilton-Ivey estimates in $3$-dim Ricci flow}
In the following computations,  we will use some cut-off functions which are composition of  a  cut-off function of $\mathbb R$ and distance function. A cut-off function $\varphi$ on real line $\mathbb R$, is a smooth nonnegative nonincreasing function, it is $1$ on $(-\infty, 1]$ and $0$ on $[2, \infty)$. We can further assume that
\begin{equation}\label{2.1.1}
{|\varphi '|\leq 2 \ , \quad |\varphi ''|+ \frac{(\varphi ')^2}{\varphi}\leq 16.
}
\end{equation}
Another often used notation is $\square= \frac{\partial}{\partial t}- \Delta$, where $\Delta$ is the Laplacian with the metric $g(t)$.

\begin{proof} of Theorem \ref{thm 1.1}. To prove the Theorem, we claim that it suffices to prove that the condition
$$
R+a(\mu+\nu)+b\nu\geq 0
$$
is preserved under complete solutions of Ricci flow for any $a,b\geq 0.$

Indeed, suppose the claim has been proved, we only need to consider the case whether $a(\mu+\nu)+b\nu\geq 0$ is preserved. If we have $a(\mu+\nu)+b\nu\geq 0$ initially, then for any $\varepsilon >0,$ we know $\varepsilon R+a(\mu+\nu)+b\nu\geq 0$ initially, hence  $\varepsilon R+a(\mu+\nu)+b\nu\geq 0$ for $t>0$ by the claim. By the arbitrariness of $\varepsilon,$  we know $a(\mu+\nu)+b\nu\geq 0$ for $t>0.$

{To prove the claim, we argue by contradiction. If there exists $(x, t)\in M^3\times (0, T]$ such that $R+a(\mu+\nu)+b\nu$ is negative at this point.  By assumption of $[R+a(\mu+\nu)+b\nu](g(0))\geq 0$ and Proposition $2.1$ in \cite{Chen}, we know $R(x, t)\geq 0$ on $[0, T]$, where $R=\lambda+\mu+\nu$ is the scalar curvature. This implies $a+b>0$. Then there are $0\leq s< s'\leq 1$ such that $s'- s\leq \frac{1}{100(a+b)}$, and the followings are satisfied:
\begin{equation}\label{2.2.1}
{\left.
\begin{array}{rl}
R+ s[a(\mu+\nu)+b\nu] &= \lambda+ (as+1)(\mu+\nu)+bs\nu\\
&= \lambda+ (a_1+1)(\mu+\nu)+b_1\nu\geq 0,
\end{array} \right.
}
\end{equation}
on $M^3\times [0, T]$, where $a_1=as$, $b_1=bs$, and

\begin{equation}\label{2.2.2}
{\left.
\begin{array}{rl}
R+ s'[a(\mu+\nu)+b\nu] &= \lambda+ (as'+1)(\mu+\nu)+bs'\nu\\
&= \lambda+ (a_1+s_1+1)(\mu+\nu)+(b_1+s_2)\nu<0,
\end{array} \right.
}
\end{equation}
at some $(x_1, t_1)$ in $M^3\times (0, T]$, where $0\leq s_1=as'-a_1\leq \frac{1}{100}$, $0\leq s_2=bs'-b_1\leq \frac{1}{100}$.

By $(\ref{2.2.1})$ and $(\ref{2.2.2})$, we know $s_1+s_2>0$, $\nu<0$ and $\lambda\geq 0$.

Choose $x_0\in M^3, r_0>0$ such that $Rc(x, t)\leq 2r_0^{-2}$ for $x\in B_{t}(x_0, r_0)$, $t\in [0, T]$. By \cite{Pere}, we have
\begin{equation}\label{2.2.3}
{\square d_t(x_0, x)\geq -\frac{10}{3}r_0^{-1},
}
\end{equation}
whenever $d_t(x_0, x)> r_0$ in the sense of support functions. Define
\begin{equation}\label{2.2.4}
{P_{ij}= \varphi\Big(\frac{d_t(x_0, x)}{Ar_0}\Big) \mathscr M_{ij},
}
\end{equation}
where $A> 1$ is any big enough number such that $(x_1, t_1)\in B_{t_1}(x_0, Ar_0),$
where $\mathscr M(x, t)_{ij}=Rg_{ij}- 2R_{ij}$ is the curvature operator matric.

In the following, we use the Uhlenbeck trick as in \cite{H4} to write the equation of the curvature operator in moving frames:
\begin{equation}\label{2.2.4'}
\frac{\partial}{\partial t}\mathscr M=\triangle \mathscr M+\mathscr M^2+\mathscr M^{\sharp}.
\end{equation}
Let $\lambda\geq \mu\geq \nu$ be eigenvalues of $\mathscr M$,  and  $\mathbb V_1, \mathbb V_2, \mathbb V_3$ be the corresponding orthonormal  eigenvectors. Then in this basis, we have
$$
      \mathscr M= \left(
       \begin{array}{ccc}
\lambda & 0 & 0 \\
0 & \mu & 0 \\
0 & 0 & \nu
       \end{array}
    \right),
$$
and
$$
      \mathscr M^2+\mathscr M^{\sharp}= \left(
       \begin{array}{ccc}
\lambda^2+\mu\nu & 0 & 0 \\
0 & \mu^2+\lambda\nu & 0 \\
0 & 0 & \nu^2+\lambda\mu
       \end{array}
    \right).
$$

Let $u(t)= \min_{x\in M^3} [\lambda + (a_1+s_1+1)(\mu+\nu)+(b_1+s_2)\nu] \varphi (x, t)$. From $(\ref{2.2.2})$, we know $u(t_1)< 0$. Assume $u(t_3)= [\lambda + (a_1+s_1+1)(\mu+\nu)+(b_1+s_2)\nu] \varphi(x_3, t_3)= \min_{t\in [0, T]} u(t)< 0$. Then $t_3\in (0, T].$

 At $(x_3,t_3)$, let $\mathbb V_1, \mathbb V_2, \mathbb V_3$ be the orthonormal eigenvectors of $\mathscr M$ with the corresponding eigenvalues $\lambda\geq \mu\geq \nu.$

To apply the maximum principle, we parallel translate $\mathbb V_1, \mathbb V_2, \mathbb V_3$  along radial geodesics emanating from $x_3$ at time $t_3.$  These local sections, denoted by $\widetilde{\mathbb V}_1, \widetilde{\mathbb V}_2, \widetilde{\mathbb V}_3$ are orthonormal,  and \begin{equation}\label{a}\frac{\partial \widetilde{\mathbb V}_i}{\partial t}(x_3,t_3)=\nabla \widetilde{\mathbb V}_i(x_3,t_3)=\triangle \widetilde{\mathbb V}_i(x_3,t_3)=0. \end{equation}

 We define a function $\tilde{u}(x,t)= [\mathscr M(\widetilde{\mathbb V}_1, \widetilde{\mathbb V}_1)+(a_1+ s_1+ 1)\mathscr M(\widetilde{\mathbb V}_2, \widetilde{\mathbb V}_2)+ (a_1+ s_1+ 1+ b_1+ s_2)\mathscr M(\widetilde{\mathbb V}_3, \widetilde{\mathbb V}_3)]\varphi(x,t)$ near $(x_3,t_3).$

It is easy to see that at any $x\in M^3$,
\begin{equation}\nonumber
{\left.
\begin{array}{rl}
&\lambda + (a_1+s_1+1)(\mu+\nu)+(b_1+s_2)\nu\\
&= \inf\{\mathscr M(\mathbb W_1, \mathbb W_1)+(a_1+ s_1+ 1)\mathscr M(\mathbb W_2, \mathbb W_2)+ (a_1+ s_1+ 1+ b_1+ s_2)\mathscr M(\mathbb W_3, \mathbb W_3)\\
& \quad \quad |\ where\ \{\mathbb W_1, \mathbb W_2, \mathbb W_3\} \ is\ an\ orthonormal\ basis\ of\ T_{x}M^3\}.
\end{array} \right.
}
\end{equation}
From this, we know \begin{equation}\label{b}\tilde{u}(x,t)\geq [\lambda + (a_1+s_1+1)(\mu+\nu)+(b_1+s_2)\nu]\varphi(x,t),\end{equation} and the equality holds at $(x_3,t_3).$

On the other hand, from (\ref{2.2.4}), (\ref{2.2.4'}) and (\ref{a}), we have
\begin{equation}\label{2.2.6}
{\left. \begin{array}{rl} \square P_{ij}& = -2 \nabla \varphi \nabla \mathscr M_{ij}+ Q_{ij},\\
\square \tilde{u}(x_3,t_3)&=-2\nabla \varphi \nabla (\frac{\tilde{u}}{\varphi})+Q(\mathbb V_1, \mathbb V_1)+(a_1+s_1+1)Q(\mathbb V_2, \mathbb V_2)\\& +(b_1+s_2)Q(\mathbb V_3, \mathbb V_3)\end{array}\right.},
\end{equation}
where
\begin{equation}\nonumber
{\left.
\begin{array}{rl}
& Q(\mathbb V_1, \mathbb V_1)+(a_1+s_1+1)Q(\mathbb V_2, \mathbb V_2)+(b_1+s_2)Q(\mathbb V_3, \mathbb V_3) \\
&\quad = \varphi [\lambda^2+ \mu\nu+ (a_1+s_1+1)[\mu^2+\nu^2+\lambda(\mu+\nu)]+(b_1+s_2)(\nu^2+\lambda\mu)\\
&\quad \quad + \frac{1}{Ar_0}[\varphi '\cdot \square d_t(x_0, x)- \varphi ''\frac{1}{Ar_0}] [\lambda + (a_1+s_1+1)(\mu+\nu)+(b_1+s_2)\nu]\\
& \quad =\varphi I+ II.
\end{array} \right.
}
\end{equation}
  At $(x_3, t_3),$ we have
\begin{equation}\nonumber
{\left.
\begin{array}{rl}
I(x_3, t_3)&=\lambda^2+ \mu\nu+\frac{(a_1+s_1+1)}{2}[\mu^2+\nu^2]+\frac{(a_1+s_1+1)}{2}[\mu^2+\nu^2]+ (a_1+s_1+1)\lambda(\mu+\nu)\\
&\quad +(b_1+s_2)\nu^2+(b_1+s_2)\lambda(\mu-\nu)+(b_1+s_2)\lambda\nu\\
&\geq \lambda[\lambda+ (a_1+1)(\mu+\nu)+b_1\nu] \\
& \quad +\frac{(a_1+s_1+1)}{2}[\mu^2+\nu^2]+(b_1+s_2)\nu^2+(b_1+s_2)\lambda(\mu-\nu) \\
& \quad +s_1\lambda(\mu+\nu)+s_2\lambda\nu.
\end{array} \right.
}
\end{equation}

Now by $(\ref{2.2.2})$, $\lambda< (a_1+s_1+1)|\mu+\nu|+(b_1+s_2)|\nu|$, so
\begin{equation}\label{2.2.10}
{\left.
\begin{array}{rl}
& |s_1\lambda(\mu+\nu)+s_2\lambda\nu|\\
&\leq s_1(a_1+s_1+1)(\mu+\nu)^2+s_2(b_1+s_2)\nu^2\\
&\quad +[s_1(b_1+s_2)+s_2(a_1+1+s_1)]|\nu||\mu+\nu|\\
&\leq (\frac{a_1+1+s_1}{50}+\frac{b_1+s_2}{200})(\mu+\nu)^2+ (\frac{a_1+1+s_1}{200}+\frac{b_1+s_2}{50})\nu^2\\
&\leq (\frac{a_1+1+s_1}{25}+\frac{b_1+s_2}{100})\mu^2+(\frac{a_1+1+s_1}{10}+\frac{b_1+s_2}{25})\nu^2.
\end{array} \right.
}
\end{equation}

Therefore
\begin{equation}\nonumber
{\left.
\begin{array}{rl}
I(x_3, t_3)&\geq \frac{a_1+1+s_1}{4}\mu^2+(\frac{a_1+1+s_1}{4}+\frac{b_1+s_2}{2})\nu^2\\
&\quad +(b_1+s_2)\lambda(\mu-\nu)-\frac{b_1+s_2}{100}\mu^2\\
&=\frac{a_1+1+s_1}{4}(\mu^2+\nu^2)+\frac{b_1+s_2}{4}\nu^2\\
&\quad +(b_1+s_2)\lambda(\mu-\nu)+\frac{b_1+s_2}{4}\nu^2-\frac{b_1+s_2}{100}\mu^2.
\end{array} \right.
}
\end{equation}

If $\mu\geq 0$, then $\lambda(\mu-\nu)\geq\mu^2$, otherwise, $\mu<0$, then $\nu^2\geq\mu^2$, so
\begin{equation}\nonumber
{\left.
\begin{array}{rl}
I(x_3, t_3)&\geq\frac{a_1+1+s_1}{4}(\mu^2+\nu^2)+\frac{b_1+s_2}{4}\nu^2\\
&\geq \frac{a_1+1+s_1}{8}(\mu+\nu)^2+\frac{b_1+s_2}{4}\nu^2\\
&\geq \frac{1}{16(a_1+1+s_1+b_1+s_2)}[(a_1+1+s_1)(\mu+\nu)+(b_1+s_2)\nu]^2\\
&\geq \frac{1}{16(a_1+b_1+2)}[(a_1+1+s_1)(\mu+\nu)+(b_1+s_2)\nu]^2\\
&\geq \frac{1}{32(a_1+b_1+2)}\{[(a_1+1+s_1)(\mu+\nu)+(b_1+s_2)\nu]^2+\lambda^2\}\\
&\geq \frac{1}{64(a_1+b_1+2)}[\lambda+(a_1+1+s_1)(\mu+\nu)+(b_1+s_2)\nu]^2.
\end{array} \right.
}
\end{equation}
Hence, at the point $(x_3,t_3),$ we have
\begin{equation}\nonumber
{\left.
\begin{array}{rl}
& Q(\mathbb V_1, \mathbb V_1)+(a_1+s_1+1)Q(\mathbb V_2, \mathbb V_2)+(b_1+s_2)Q(\mathbb V_3, \mathbb V_3) \\
&\quad \geq \frac{1}{64(a_1+b_1+2)\varphi} \Big [u^2- \frac{64(a_1+b_1+2)}{Ar_0}(\frac{10\varphi '}{3r_0}+ \frac{\varphi ''}{Ar_0})u \Big].
\end{array} \right.
}
\end{equation}
Combining (\ref{b}), (\ref{2.2.6}) with the fact $u(t_3)< 0$, and applying the maximum principle, we have

\begin{equation}\nonumber
{\left.
\begin{array}{rl}
0&\geq \frac{1}{64(a_1+b_1+2)\varphi} \Big [u^2- \frac{64(a_1+b_1+2)}{Ar_0}(\frac{10\varphi '}{3r_0}+ \frac{\varphi ''}{Ar_0})u \Big]+ \frac{2}{(Ar_0)^2}(\frac{\varphi '}{\varphi})^2 u(t_3)  \\
& \geq \frac{1}{64(a_1+b_1+2)\varphi} \Big [u^2+ \frac{128(a_1+b_1+2)}{Ar_0}(|\frac{10\varphi '}{3r_0}|+ |\frac{\varphi ''}{Ar_0}+ \frac{\varphi '^2}{Ar_0 \varphi}|) u(t_3)  \Big].
\end{array} \right.
}
\end{equation}
If $u(t_3)< -\frac{256(a_1+b_1+2)}{Ar_0} (\frac{20}{3r_0}+ \frac{16}{Ar_0})$, we get
\begin{equation}\label{2.2.13}
0 \geq \frac{1}{128(a_1+b_1+2)} u(t_3)^2> 0
\end{equation}
 which is a contradiction. Hence $u(t_3)\geq -\frac{256(a_1+b_1+2)}{Ar_0} (\frac{20}{3r_0}+ \frac{16}{Ar_0})$ on $[0, T]$. By the definition of $u(t)$, we get
\begin{equation}\label{2.2.14}
{\lambda+ (a_1+s_1+1)(\mu+\nu)+(b_1+s_2)\nu\geq -\frac{1}{Ar_0}C(a_1+b_1, r_0)
}
\end{equation}
on $B_t(x_0, Ar_0)$, for $t\in [0, T]$. Then let $A\rightarrow \infty$, we get $\lambda+ (a_1+s_1+1)(\mu+\nu)+(b_1+s_2)\nu\geq 0$ on $M^3\times [0, T]$. That is contradiction with (\ref{2.2.2}). The proof of Theorem \ref{thm 1.1} is completed.
}

\end{proof}

By taking $a=0,b=1,c=0$ in Theorem \ref{thm 1.1}, we get Theorem \ref{thm 1.2}.

\begin{remark}\label{remark 2.3}
{In $3$-dim Ricci flow, nonnegative Ricci curvature is always preserved with bounded curvature assumption, which can be proved by using the maximum principle directly. The above theorem removes the bounded curvature assumption. In \cite{Chen}, the local estimate of $Rm$ is achieved by an induction method. Here we make the argument in a `continuous' way.
}
\end{remark}

In the following,  we will prove the local Hamilton-Ivey estimate without bounded curvature assumption on  $3$-dim manifolds.

\begin{proof} of Theorem \ref{thm 1.3}.
{Assume that $\lambda\geq \mu\geq \nu$ are the eigenvalues of curvature operator $\mathscr M_{ij}$ and  $R= \lambda+ \mu+ \nu$ as before. Note our curvature operator's value on tangent plane is two times sectional curvature of the tangent plane.

% By contradiction. If there exists $(x, t)\in M^3\times [0, T]$, (\ref{2.3.1}) does not hold, it is easy to get $t> 0$.

Let $w= \frac{R}{-\nu}- \ln (-\nu)+ 3+ \ln\frac{K}{1+ Kt}$ be a function defined on
\[\Omega= \{(x, t)|\ \nu(x, t)< 0, \ (x, t)\in M^3\times [0, T]\}\]
We want to prove (\ref{1.3.1}) on $\Omega\cap B_{t}(x_0, Ar_0)$.  Recall that as in (\ref{2.2.4'}), we  pick a abstract vector bundle $E$ with fixed bundle metric and a family of time-dependent connections compatible with the metric.   After pulling back with the moving frames from the tangent bundle, the curvature operator $\mathscr M$ acts on $E$.

Define  \begin{equation}\label{2.3.1.2.1}
{\left.
\begin{array}{rl}
W(x, t, \mathbb{V})&= -R\cdot [\mathscr M(x, t)(\mathbb{V}, \mathbb{V})]^{-1}- \ln (-\mathscr M(x, t)(\mathbb{V}, \mathbb{V}))\\ &\quad + \Big( 3+ \ln\frac{K}{1+ Kt} \Big)\cdot |\mathbb{V}|^2
\end{array} \right.
}
\end{equation}
to be a function on
\[\Omega_1= \{(x, t, \mathbb{V})|\ \mathscr M(x, t)(\mathbb{V}, \mathbb{V})< 0,\ (x, t)\in \Omega,\ \mathbb{V}\in E, \ |\mathbb{V}|= 1\}.\]

When $R> \nu$, considering the function $h_1(s)= \frac{R}{-s}- \ln(-s)$, we find it is increasing in $s$. Hence, when $R> \nu$, $(x, t)\in \Omega$, we have
\begin{equation}\label{2.3.1.3}
w(x, t) = \min \left\{
\begin{array}{c}
W(x, t,\mathbb V) ;\\
 where\ (x, t, \mathbb{V})\in \Omega_1
\end{array} \right\}.
\end{equation}

Note $R(x, 0)\geq 3\nu(x, 0)\geq -3K$ on $B_0(x_0, 2Ar_0)$, by Proposition $2.1$ of \cite{Chen}, we have
\begin{equation}\label{2.3.5.1}
{R(x, t)\geq \min\Big\{-\frac{3K}{1+ Kt}, -\frac{C}{Ar_0^2}\Big\},
}
\end{equation}
where $x\in B_{t}(x_0, \frac{3}{2}Ar_0)$, $t\in [0, T]$, and $C$ is some universal constant. By the assumption $A\geq \frac{C}{3r_0^2}\Big(\frac{1}{K}+ T\Big)$ and (\ref{2.3.5.1}),  we  get  $$R\geq -\frac{3K}{1+ Kt},$$ on $B_{t}(x_0, Ar_0)$ .

If $0> \nu\geq R$, assume $\frac{R}{-\nu}= a$, then $-3\leq a\leq -1$. Note $-\nu= \frac{1}{a}R\leq - \frac{3}{a}\frac{K}{1+ Kt}$, hence
\begin{equation}\label{2.3.1.2}
{\left.
\begin{array}{rl}
w\geq a+ \ln (-a)+ 3- \ln 3\geq 0.
\end{array} \right.
}
\end{equation}
The last inequality uses the fact that the minimum of function $a+ \ln (-a)$ on $[-3, -1]$ is $(-3+ \ln 3)$.
By (\ref{2.3.1.2}), we only need to show (\ref{1.3.1}) for the case $R> \nu$.

Let
\begin{equation}\label{2.3.2}
{u(x, t)= \varphi\Big(\frac{2d_t(x_0, x)}{Ar_0}- 1\Big)w(x, t)
}
\end{equation}
be a function defined on $\Omega_1$. Let $(x_3, t_3, \mathbb{V}_3)$ be a point such that
\begin{equation}\label{2.3.5}
\varphi(x_3,t_3)W(x_3,t_3,\mathbb{V}_3)= \min_{t\in [0, T],x\in \Omega}(\varphi w)(x,t).
\end{equation}
If $\varphi(x_3,t_3)W(x_3,t_3,\mathbb{V}_3)\geq 0$, we are done.

 So we may assume $\varphi(x_3,t_3)W(x_3,t_3,\mathbb{V}_3)< 0$, then $t_3\in (0, T]$ and $x_3\in B_{t_3}(x_0, \frac{3}{2}Ar_0)$.  Note that we have $R(x_3, t_3)\geq -\frac{3K}{1+ Kt}$.

At $(x_3, t_3, \mathbb{V}_3)$, let $\widetilde{\mathbb{V}}$  be a local vector field defined by parallel translations of $\mathbb{V}_3$ along radial geodesics emanating from $x_3$ at time $t_3.$ Then $|\widetilde{\mathbb{V}}|\equiv1$ and  $\frac{\partial}{\partial t}\widetilde{\mathbb{V}}=\nabla \widetilde{\mathbb{V}}=\triangle \widetilde{\mathbb{V}}=0$ at $(x_3,t_3).$ Define two smooth functions $$\tilde{\nu}=\mathscr M(x, t)(\widetilde{\mathbb{V}}, \widetilde{\mathbb{V}})$$ $$\tilde{u}=\Big[\frac{R}{-\tilde{\nu}}- \ln (-\tilde{\nu})+ 3+ \ln\frac{K}{1+ Kt} \Big]\varphi(x, t)$$ near $(x_3,t_3),$ these two functions satisfy
\begin{equation}\label{2.3.5'}
\begin{array}{c}
\tilde{\nu}\geq \nu, \tilde{u}\geq u
\end{array}
\end{equation}
and equalities hold at $(x_3,t_3).$

At $(x_3,t_3),$ a straightforward computation gives
\begin{equation}\label{2.3.4}
{\square \tilde{u}= \Big[2\frac{\varphi '}{Ar_0}\square d_t- 4\frac{\varphi ''}{(Ar_0)^2}\Big]W(\cdot,\cdot,\widetilde{\mathbb{V}})- 2\nabla \varphi \nabla W(\cdot,\cdot, \widetilde{\mathbb{V}})+ \varphi \square W(\cdot,\cdot, \widetilde{\mathbb{V}})
}
\end{equation}and
\begin{equation}\label{2.3.3}
{\left.
\begin{array}{rl}
 \square W(x, t, \widetilde{\mathbb{V}})&= 2\nabla \ln(-\tilde{\nu})\cdot \nabla \tilde{w}+ |\nabla \ln(-\tilde{\nu})|^2 \\
&\quad - 2\tilde{\nu}- \frac{K}{1+ Kt}+ \nu^{-2}(\lambda- \nu)(\mu- \nu)R.
\end{array} \right.
}
\end{equation}

 From(\ref{2.3.4}) and (\ref{2.3.3}), we have
\begin{equation}\label{2.3.7}
{\left.
\begin{array}{rl}
\square \tilde{u}&\geq \varphi\Big[-2\nu-\frac{K}{1+ Kt}\Big] + \varphi \frac{(\lambda- \nu)(\mu- \nu)}{\nu^2}R\\
&\quad +4\frac{2}{Ar_0^2}\Big(|\frac{10\varphi '}{3}|+ |\frac{\varphi ''}{A}|+ |\frac{1}{A}\cdot \frac{(\varphi ')^2}{\varphi}|\Big)w- \frac{64}{(Ar_0)^2}w^2.
\end{array} \right.
}
\end{equation}
Combining with (\ref{2.3.5}) and (\ref{2.3.5'}), and applying the maximum principle at $(x_3, t_3),$  we get
\begin{equation}\label{2.3.7.1}
{\left.
\begin{array}{rl}
0& \geq \square \tilde{u}|_{(x_3, t_3)}\\
&\geq \varphi\Big[-2\nu-\frac{K}{1+ Kt}\Big] + \varphi\frac{(\lambda- \nu)(\mu- \nu)}{\nu^2}R+ \frac{120}{Ar_0^2}w- \frac{64}{(Ar_0)^2}w^2\\
&= (I)+ (II)+ (III)+ (IV).
\end{array} \right.
}
\end{equation}

The rest computations are all at $(x_3, t_3)$ without illustration.We estimate $(II)$ firstly. If $R(x_3, t_3)\geq 0$, then $(II)\geq 0$.

If $R(x_3, t_3)< 0$, we have
\begin{equation}\label{2.3.8}
{\left.
\begin{array}{rl}
(II)&\geq \varphi\frac{(R- 3\nu)^2}{4\nu^2}R\geq \frac{\varphi}{4}\Big(\frac{R}{-\nu}+ 3\Big)^2 R\\
&\geq \frac{\varphi}{4}\Big[\ln(-\nu)- \ln(\frac{K}{1+ Kt})\Big]^2 R.
\end{array} \right.
}
\end{equation}

%Next we estimate $(III)$, using the property of $\varphi$, we get
%\begin{equation}\label{2.3.9}
%{(III)\geq \frac{C(r_0)}{A}\varphi^{-1}u(t_3)
%}
%\end{equation}

Now we will estimate $(I)+ (II)$ in terms of $\nu$. Note
\begin{equation}\label{2.3.9}
{0> w= \frac{R}{-\nu}- \ln(-\nu)+ 3+ \ln\frac{K}{1+ Kt_3}\geq -\ln(-\nu)+ \ln\frac{K}{1+ Kt_3}
},
\end{equation}
hence we get
\begin{equation}\label{2.3.9.1}
{|w(x_3, t_3)|\leq \Big|\ln (-\nu)- \ln\frac{K}{1+ Kt_3}\Big|.
}
\end{equation}
On the other side, by $w(x_3, t_3)< 0$, we get
\begin{equation}\label{2.3.10}
{\ln (-\nu)> \frac{R}{-\nu}+ 3+ \ln\frac{K}{1+ Kt_3}\geq \ln\frac{K}{1+ Kt_3}.
}
\end{equation}
Then $-\nu> \frac{K}{1+ Kt_3}$, and we get
\begin{equation}\label{2.3.10.1}
{(I)\geq \varphi \cdot(-\nu).
}
\end{equation}
Now consider the function
\begin{equation}\nonumber
{f(s)= \frac{s}{\Big|\ln s- \ln\frac{K}{1+ Kt_3}\Big|^2}
}
\end{equation}
where $s> \frac{K}{1+ Kt_3}$. Then it is easy to get $f(s)\geq f(e^2\cdot \frac{K}{1+ Kt_3})= \frac{e^2}{4}\cdot \frac{K}{1+ Kt_3}$ for any $s> \frac{K}{1+ Kt_3}$. From this and (\ref{2.3.9.1}), we get
\begin{equation}\label{2.3.11}
{-\nu\geq \frac{e^2}{4}\cdot\frac{K}{1+ Kt_3}\Big[\ln(-\nu)- \ln\frac{K}{1+ Kt_3}\Big]^2\geq \frac{e^2}{4}\cdot\frac{K}{1+ Kt_3}w^2.
}
\end{equation}
Hence
\begin{equation}\label{2.3.11.1}
{(II)\geq \frac{\varphi}{e^2}(-\nu)\cdot\frac{1+ Kt_3}{K}R\geq \frac{3}{e^2}\varphi \nu.
}
\end{equation}
     As a consequence,   we always have
\begin{equation}\label{2.3.12}
{(I)+ (II)\geq -\frac{1}{9}\varphi \nu.
}
\end{equation}
To control $(IV)$, we consider the function
\begin{equation}\nonumber
{h_2(s)= \frac{s}{\Big|\ln s- \ln\frac{K}{1+ Kt_3}\Big|^3}
}
\end{equation}
where $s> \frac{K}{1+ Kt_3}$. Then it is easy to get $h_2(s)\geq h_2(\frac{Ke^3}{1+ Kt_3})= \frac{e^3}{27}\cdot\frac{K}{1+ Kt_3}$ for any $s> \frac{K}{1+ Kt_3}$. As in (\ref{2.3.11}), we get
\begin{equation}\label{2.3.12.1}
{-\nu\geq \frac{e^3}{27}\frac{K}{1+ Kt_3} |w|^3.
}
\end{equation}
By (\ref{2.3.11}), (\ref{2.3.12}) and (\ref{2.3.12.1}), we have
\begin{equation}\label{2.3.12.2}
{(I)+ (II)\geq -\frac{1}{9}\varphi \nu\geq \frac{e^2}{72}\cdot\frac{K}{1+ Kt_3} \varphi w^2+ \frac{e^3}{486}\cdot\frac{K}{1+ Kt_3} \varphi |w|^3.
}
\end{equation}

%If $R(x_3, t_3)\geq 0$, by (\ref{2.3.7.1}), (\ref{2.3.12}), we get
%\begin{equation}\label{2.3.13}
%{\left.
%\begin{array}{rl}
%0\geq \frac{e^2}{4\varphi}\cdot\frac{K}{1+ Kt_3}V^2+ \frac{120}{Ar_0\varphi}V
%\end{array} \right.
%}
%\end{equation}
%By $u(x_3, t_3)< 0$, we get $u(x_3, t_3)\geq -\frac{100(1+ Kt_3)}{AKr_0}$.

%If $R(x_3, t_3)< 0$, by (\ref{2.3.7.1}), (\ref{2.3.8}) and (\ref{2.3.12}), we have

By (\ref{2.3.7.1}) and (\ref{2.3.12.2}), we get
\begin{equation}\label{2.3.14}
{\left.
\begin{array}{rl}
0&\geq \frac{e^2}{72\varphi}\cdot\frac{K}{1+ Kt_3}(\varphi w)^2+ \frac{120}{Ar_0^2\varphi}\varphi w \\
&\quad + \frac{e^3}{486\varphi}\cdot\frac{K}{1+ Kt_3}|w|\Big( (\varphi w)^2+ \frac{486}{e^3}\cdot\frac{1+ Kt_3}{ K}\cdot \frac{64}{(Ar_0)^2}\varphi w \Big)
\end{array} \right.
}
\end{equation}

% If $|V|\geq \frac{960(1+ Kt_3)}{e^2KAr_0}$, we get
%\begin{equation}\label{2.3.15}
%{\left.
%\begin{array}{rl}
%0\geq \frac{e^2}{8\varphi}\cdot\frac{K}{1+ Kt_3}V^2+ \frac{9}{4}\varphi R
%\end{array} \right.
%}
%\end{equation}
%Hence by (\ref{2.3.5.1}), we have
%\begin{equation}\label{2.3.16}
%{\left.
%\begin{array}{rl}
%V^2\leq \frac{18(1+ Kt_3)}{e^2K}\cdot(-R)\leq \max\Big\{\frac{54}{e^2},\ \frac{C(1+ Kt_3)}{KAr_0^2}\Big\}
%\end{array} \right.
%}
%\end{equation}

From  the above, we get
\begin{equation}\label{2.3.17}
{\left.
\begin{array}{rl}
|\varphi w|\leq \max\Big\{\frac{8640(1+ Kt_3)}{KAe^2r_0^2},\ \frac{31104(1+ Kt_3)}{Ke^3(Ar_0)^2} \Big\} \\
\end{array} \right.
}
\end{equation}

By $\varphi= 1$ on $B_{t}(x_0, Ar_0)$, $t\in [0, T]$ and (\ref{2.3.17}), we have
\begin{equation}\label{2.3.18}
{\left.
\begin{array}{rl}
\frac{R}{-\nu}- \ln(-\nu)+ 3+ \ln\frac{K}{1+ Kt}\geq \min\Big\{-\frac{8640(1+ KT)}{KAe^2r_0^2},\ -\frac{31104(1+ KT)}{Ke^3(Ar_0)^2} \Big\}.
\end{array} \right.
}
\end{equation}

Theorem \ref{thm 1.3} is proved.
}
\end{proof}

\begin{remark}\label{remark 2.5}{ Our method of using maximal principle here is a little bit different from Hamilton's way of using maximal principle in \cite{H4} and \cite{Ha99}. The estimate similar to (\ref{2.3.12}) had appeared in \cite{CLN} in the context of compact manifold (see $(6.39)$ there).
}
\end{remark}
Now we give the proof of Corollary \ref{thm 1.4}.

\begin{proof} of Corollary \ref{thm 1.4}.
{We choose $x_0\in M^3$, $r_0> 0$ such that $Rc(x, t)\leq 2r_0^{-2}$ for $x\in B_{t}(x_0, r_0)$, $t\in [0, T]$. This is always possible if we pick small $r_0.$ For any constant  $A\geq \frac{C}{3r_0^2}(1+ T)$,  by Theorem \ref{thm 1.3},  we get
\begin{equation}\label{2.4.2}
{\left.
\begin{array}{rl}
\frac{R}{-\nu}- \ln(-\nu)+ 3+ \ln(\frac{K}{1+Kt})\geq \min \Big\{-\frac{8640(1+ KT)}{KAe^2r_0^2},\ -\frac{31104(1+ KT)}{Ke^3(Ar_0)^2} \Big\}.
\end{array} \right.
}
\end{equation}
at any points $(x_1, t_1)\in B_{t}(x_0,Ar_0)$ with $\nu(x_1, t_1)< 0.$

Let $A\rightarrow \infty$ in (\ref{2.4.2}), we get
\begin{equation}\label{2.4.3}
{\left.
\begin{array}{rl}
\frac{R}{-\nu}- \ln(-\nu)+ 3+  \ln(\frac{K}{1+Kt})\geq 0.
\end{array} \right.
}
\end{equation}
}
\end{proof}

If we have a ancient solution defined on $(-\infty, 0],$  and there were some point $(x_1,t_1)$ such that $\nu(x_1, t_1)< 0,$ then by Corollary \ref{thm 1.5}, we have \begin{equation}\nonumber
R(x_1,t_1)\geq (-\nu)(x_1,t_1)(-3+\ln (t_1-t_2)(-\nu)(x_1,t_1)))
\end{equation}for any $t_2<t_1$. When $t_2\rightarrow -\infty,$ this will give a contradiction. Therefore, we have the following pinching result about ancient solution originally due to the first author (see \cite{Chen}).

\begin{cor}\label{cor 2.5}
{Any smooth complete ancient solution of Ricci flow on $3$-dim manifold must have $Rm\geq 0$.
}
\end{cor}

Following similar strategy in \cite{Chen}, we also have the following theorem.
\begin{theorem}\label{thm 4.1}
{Let $(M^3, g(x))$ be a complete noncompact $3$-dim manifold with $0\leq Rc\leq K_0g$, for some fixed positive constant $K_0$. Also assume $i_0= \min_{x\in M^3} inj(x)\geq \delta> 0$. If $g_1(t)$ and $g_2(t)$ are both smooth complete solutions to the Ricci flow on $M^3\times [0, T]$ with $g(x)$ as initial data, we have $g_1(t)\equiv g_2(t)$, for $0\leq t< \min\{T, \frac{1}{45K_0}\}$.
}
\end{theorem}

\section*{Acknowledgements}
B.-L. Chen was partially supported by NSFC 11025107, and High Level Talent Project in High Schools in Gongdong Province (34000-5221001). G. Xu thanks Bo Yang in University of California at San Diego for many interesting discussions. Z. Zhang was partially supported by NSFC 11126204, foundation for Distinguished Young Talents in Higher Education of Guangdong (2012LYM-0051) and Youth Foundation of South China Normal University.

\end{document}